\documentclass[10pt,twocolumn,twoside]{IEEEtran}

\usepackage[margin=54pt]{geometry}

\usepackage{mathrsfs}
\usepackage{amsfonts}
\usepackage{amssymb}
\usepackage{amsmath}
\usepackage{graphicx}
\usepackage{cite}
\usepackage{amsthm}
\usepackage{commath}
\usepackage{tikz}
\usepackage{tkz-graph}
\usepackage[europeanresistors]{circuitikz}
\usetikzlibrary{arrows}
\usepackage{psfrag}
\usepackage{pgfplots}
\usetikzlibrary{external}
\usepackage{booktabs}
\usepackage{enumerate}
\allowdisplaybreaks
\usepackage{float}
\usepackage[lofdepth,lotdepth]{subfig}
\usepackage[font = small]{caption}
\pgfplotscreateplotcyclelist{Power_frequency}{%
	{blue,line width=1pt},
	{red,line width=1pt},
	{magenta,line width=1pt},
	{yellow,line width=1pt},		
	{cyan,line width=1pt},
}

\newtheorem{theorem}{Theorem}
\newtheorem{corollary}[theorem]{Corollary}
\newtheorem{lemma}[theorem]{Lemma}

\newtheorem{definition}{Definition}
\newtheorem{assumption}{Assumption}

\DeclareMathOperator*{\diag}{diag}

\DeclareMathOperator*{\E}{E}

\DeclareMathOperator*{\tr}{tr}

\newcommand{\beq}{\begin{equation}}
\newcommand{\eeq}{\end{equation}}
\newcommand{\bq}{\begin{eqnarray}}
\newcommand{\eq}{\end{eqnarray}}
\newcommand{\bqn}{\begin{eqnarray*}}
\newcommand{\eqn}{\end{eqnarray*}}
\newcommand{\bee}{\begin{enumerate}}
\newcommand{\eee}{\end{enumerate}}

\newcommand{\hn}{$\mathcal{H}_2$ }

\newboolean{showcomments}
\setboolean{showcomments}{true}

\newcommand{\martin}[1]{\ifthenelse{\boolean{showcomments}}
{\textcolor{blue}{Martin says: #1}}{}}
\newcommand{\emma}[1]{\ifthenelse{\boolean{showcomments}}
{\textcolor{magenta}{Emma says: #1}}{}}

\newlength\fheight
\newlength\fwidth

\begin{document}
\title{\vspace*{18pt} Performance and Scalability of Voltage Controllers in Multi-Terminal HVDC Networks }

\author{ \IEEEauthorblockA{Martin Andreasson, Emma Tegling, Henrik Sandberg and Karl H. Johansson}
 
 \thanks{The authors are with the School of Electrical Engineering and the ACCESS Linnaeus Centre, KTH Royal Institute of Technology, SE-100 44 Stockholm, Sweden {(\tt mandreas, tegling, hsan, kallej@kth.se)}.}
\thanks{This work was supported in part by the European Commission, the Knut and Alice Wallenberg Foundation, by the Swedish Research Council through grants 2014-6282 and 2013-5523, and the Swedish Foundation for Strategic Research.}
}


        
\maketitle
\thispagestyle{empty}
\pagestyle{empty}

\begingroup
\makeatletter
\renewcommand{\p@subfigure}{}
\makeatother

\begin{abstract}
In this paper, we compare the transient performance of a multi-terminal high-voltage DC (MTDC) grid equipped with a slack bus for voltage control to that of two distributed control schemes: a standard droop controller and a distributed averaging proportional-integral (DAPI) controller. We evaluate performance in terms of an \hn metric that quantifies expected deviations from nominal voltages, and show that the transient performance of a droop or DAPI controlled MTDC grid is always superior to that of an MTDC grid with a slack bus. In particular, by studying systems built up over lattice networks,  we show that the \hn norm of a slack bus controlled system may scale unboundedly with network size, while the norm remains uniformly bounded with droop or DAPI control. 
 We simulate the control strategies on radial MTDC networks to demonstrate that the transient performance for the slack bus controlled system deteriorates significantly as the network grows, which is not the case with the distributed control strategies. 
\end{abstract}

\section{Introduction}
Transmitting power over long distances while maintaining low losses is one of the greatest challenges related to power transmission systems. Driven partly by increased deployment of renewable energy resources, such as large-scale off-shore wind farms, many distances between power generation and consumption are increasing. There is therefore a growing need for long-distance power transmission, motivating a widespread use of high-voltage direct current (HVDC) technology. Its higher investment costs compared to AC transmission lines are compensated by its lower resistive losses for sufficiently long distances, which are typically 500-800 km for overhead lines \cite{padiyar1990hvdc}, but less than 100 km for undersea cable connections \cite{bresesti2007hvdc}. 
As more energy sources and consumers are connected by HVDC lines, the individual lines will eventually form a grid consisting of multiple terminals connected by several HVDC transmission lines, resulting in so-called multi-terminal HVDC (MTDC) systems \cite{van2010multi}. 

The operation of MTDC transmission systems relies on the ability to control the DC voltages at the terminals; firstly, in order to govern the network's current flows, and secondly, in order to avoid damage to power electronic equipment caused by too large deviations from nominal operating voltages \cite{kundur1994power, van2010multi}. Different schemes for this voltage control in HVDC systems have been proposed in the literature. One method is to assign one of the buses (the \textit{slack bus}) to control the networks' voltage drift through, for example, a proportional-integral controller \cite{Haileselassie2012, Wang2014}. Remaining buses control their injected currents according to Ohm's and Kirchhoff's laws~\cite{kundur1994power}. We refer to this control strategy as \textit{slack bus control}. The well-known \textit{voltage droop controller}, on the other hand, is a decentralized proportional controller that regulates current injections based on local voltages \cite{kundur1994power, KarlssonSvensson2003, Haileselassie2012}. This, however, typically leads to stationary voltage errors that can be eliminated through so-called secondary control. Several secondary controllers have been proposed in recent work \cite{ Sarlette20123128, Anand2013, Xiaonan2014,  andreasson2014TCNS}. Here, we focus on a \textit{distributed averaging proportional-integral (DAPI)} controller. 

The control problem aspects outlined above are relevant also for DC microgrids, which have attracted research interest in recent years. DC microgrids are thought of as low-voltage distribution networks with distributed generation sources, storage elements and loads, all operating on DC~\cite{ZhaoDorfler2015,DCmicrogridsurvey }. Although this paper will focus on HVDC systems, the same analysis can in principle be applied to DC microgrids, after a network reduction procedure laid out in~\cite{ZhaoDorfler2015}. 

The objective of this paper is to analyze the transient performance of the MTDC grid and to compare the slack bus control strategy to the voltage droop and the DAPI controllers. We  evaluate performance through an \hn norm metric that quantifies each node's expected voltage deviations over the voltage regulation transient. We show that the performance of a droop or DAPI controlled MTDC grid is always superior to that of a slack bus controlled grid. We also derive theoretical bounds on the scaling of the \hn norms with network size by studying systems built up on large-scale lattice networks. We find that, while the system's \hn norm remains bounded with droop or DAPI control for any network structure, that of a slack bus controlled system grows unboundedly with network size in 1- and 2-dimensional lattice networks. Our results therefore indicate that the slack bus control strategy is \textit{not scalable} to larger networks, while droop and DAPI control are. 

The remainder of this paper is organized as follows. We introduce the MTDC network model, the voltage controllers, and the \hn performance metric in Section~\ref{sec:model}. In Section~\ref{sec:computing_H2}, we calculate the systems' \hn norms and discuss their scaling with network size in Section~\ref{sec:large-scale}. We present a numerical simulation in Section~\ref{sec:simulations} and conclude in Section~\ref{sec:discussion}. 




\section{Model and problem setup}
\label{sec:model}

\subsection{Notation}
\label{subsec:prel}
Let $\mathcal{G} = (\mathcal{V}, \mathcal{E})$ be a graph, where $\mathcal{V}=\{ 1,\hdots, n \}$ denotes the vertex set and $\mathcal{E}=\{ 1,\hdots, m \}$ denotes the edge set. Let $\mathcal{N}_i$ be the neighbor set of vertex $i \in \mathcal{V}$ in $\mathcal{G}$.
For the MTDC network considered here, $\mathcal{V}$ corresponds to the HVDC bus set, and $\mathcal{E}$ corresponds to the set of HVDC lines. 
Throughout this paper, we assume that the graph that underlies the MTDC network is connected. 
Denote by $\mathcal{B}$ the vertex-edge adjacency matrix of a graph, and let $\mathcal{L}_W=\mathcal{B}W\mathcal{B}^T$ be its weighted Laplacian matrix, with edge-weights given by the elements of the positive definite diagonal matrix $W$.
Let $c_{n\times m}$ be a matrix of dimension $n\times m$ whose elements are all equal to the number~$c$, and $c_n$ a column vector whose elements are all equal to $c$. 
Let $J_n=\frac{1}{n} 1_{n\times n}$, and denote by $A^*$ the conjugate transpose of the matrix $A$. 
For simplicity, we will often drop the time dependence of variables in the notation, e.g., $x(t)$ will be denoted $x$.

\subsection{Model}
Consider an MTDC transmission system consisting of $n$ HVDC terminals, denoted by the vertex set $ \mathcal{V}= \{1, \dots, n\}$. The DC terminals are connected by $m$ HVDC transmission lines, denoted by the edge set $ \mathcal{E}= \{1, \dots, m\}$. 
The HVDC lines are assumed to be purely resistive, neglecting any capacitive and inductive line elements. The assumption of purely resistive lines is not restrictive for the control applications considered in this paper, since line capacitance can be included in the capacitances of the terminals \cite{kundur1994power}. This implies that the current $I_{ij}$ on the HVDC line from terminal~$i$ to terminal~$j$ is given by Ohm's law as 
\begin{align*}
I_{ij} = \frac{1}{R_{ij}} (V_i -V_j),
\end{align*}
where $V_i$ is the voltage deviation from the nominal voltage $V^\text{nom}_i$ of terminal $i$, and $R_{ij}$ is the line resistance. 
The voltage dynamics of an arbitrary DC terminal~$i$ are assumed to be given by
\begin{align}
C_i \dot{V}_i &= -\sum_{j\in \mathcal{N}_i} I_{ij} + u_i = -\sum_{j\in \mathcal{N}_i} \frac{1}{R_{ij}}(V_i -V_j) +  u_i,
\label{eq:hvdc_scalar}
\end{align}
where $C_i>0$ is the total capacitance of terminal $i$, including that of the incident HVDC line as well as any shunt capacitances, and $u_i$ is the controlled injected current for which we will shortly introduce control schemes. 
Equation \eqref{eq:hvdc_scalar} may be written in vector-form as
\begin{align}
\begin{aligned}
C \dot{V} &= -\mathcal{L}_R V + u,
\end{aligned}
\label{eq:hvdc_vector}
\tag{$\mathcal{S}_\text{MTDC}$}
\end{align}
where $V=[V_1, \dots, V_n]^T$, $C=\diag([C_1, \dots, C_n])$, $u=[u_1, \dots, u_n]^T$ and $\mathcal{L}_R$ is the weighted Laplacian matrix of the MTDC network graph, with edge-weights given by the conductances $\frac{1}{R_{ij}}$. Fig.~\ref{fig:graph} illustrates a four terminal MTDC system. 

\subsection{Slack bus control}
A common control strategy for MTDC grids is to control the voltage at one terminal, by means of, e.g., a proportional-integral controller. This terminal, which then regulates the network's voltage drift, is called a slack bus. Such a control strategy can be idealized by assuming that the slack bus is grounded, that is, 
\begin{align}
\label{eq:slack_bus}
V_1(t) = 0 ,~~\forall t \ge 0,
\end{align}
where, without loss of generality, we have assigned terminal~1 to be the slack bus. This results in the following dynamics for the remaining buses
\begin{align}
\begin{aligned}
\tilde{C} \dot{\tilde{V}} &= -\tilde{\mathcal{L}}_R \tilde{V},
\end{aligned}
\label{eq:hvdc_vector_slack}
\tag{$\mathcal{S}_\text{slack}$}
\end{align}
where $\tilde{\mathcal{L}}_R$ is the reduced Laplacian matrix of ${\mathcal{L}}_R$, which is obtained by deleting the first row and the first column of ${\mathcal{L}}_R$, $\tilde{C}$ is obtained mutatis mutandis, and $\tilde{V}=[V_2, \dots, V_n]^T$.  

With slack bus control (also called constant DC voltage control), the network's operation relies on the functioning of one single terminal. Therefore, to increase the reliability of multi-terminal networks, distributed or decentralized control schemes inspired by frequency control in AC networks have been proposed \cite{Haileselassie2012, Wang2014}. We next introduce two such controllers. 

\begin{figure}
\center
\includegraphics[scale=0.9]{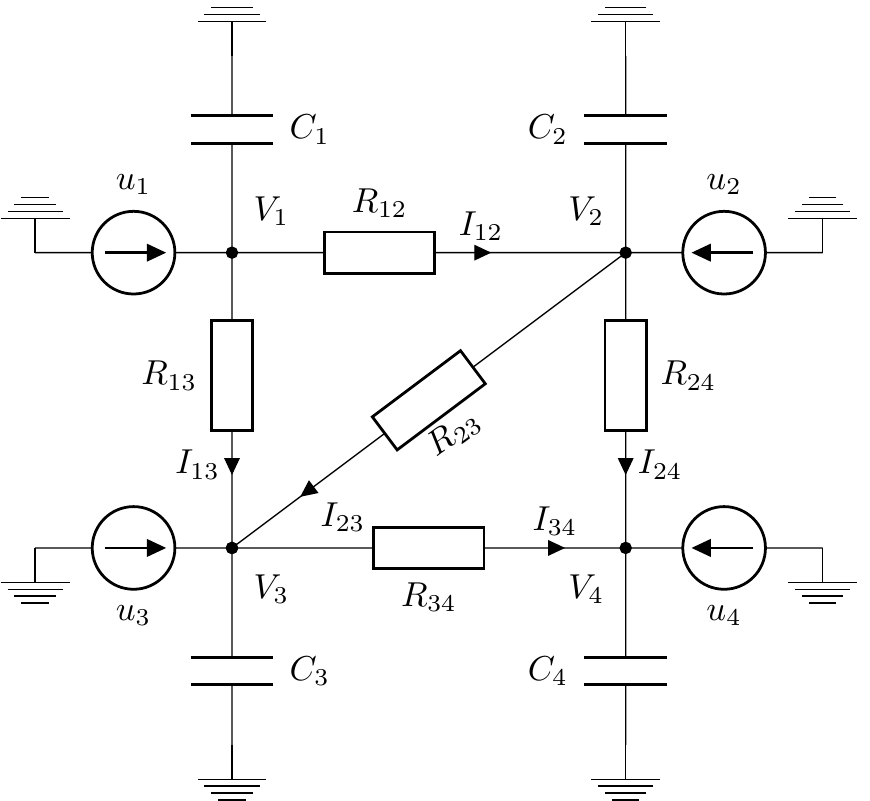}
\caption{Example of an MTDC network consisting of 4 terminals (buses) and 5 lines.}
\label{fig:graph}
\end{figure}

\subsection{Droop control}
The voltage droop controller is a commonly proposed method for controlling \eqref{eq:hvdc_vector}, see, e.g., \cite{ZhaoDorfler2015,Haileselassie2012}. It is a decentralized proportional controller, which takes the form
\begin{align}
\label{eq:droop_vector}
u &= -K_PV,
\end{align}
where $K_P = \diag\{K_{P_1}, \ldots, K_{P_n}\}$ contains the droop gains $K_{P_i}>0$. Inserting \eqref{eq:droop_vector} in \eqref{eq:hvdc_vector} yields 
\begin{align}
\begin{aligned}
C \dot{V} &= -(\mathcal{L}_R + K_P) V .
\end{aligned}
\label{eq:hvdc_vector_droop}
\tag{$\mathcal{S}_\text{droop}$}
\end{align}

\subsection{Distributed averaging proportional-integral control}
Various secondary controllers have been proposed for MTDC grids and DC microgrids, with the objective to achieve current sharing and to eliminate static voltage errors \cite{ Sarlette20123128, Anand2013, Xiaonan2014,  andreasson2014TCNS}. Here, we consider a distributed averaging proportional-integral (DAPI) controller, which appends a secondary controller layer with an associated communication network to the droop controller \eqref{eq:hvdc_vector_droop}. The DAPI controller has been successfully applied in frequency control of AC grids \cite{SimpsonPorco2013synchronization, Andreasson2013_ecc, Simpson2015Secondary, Trip2016Internal} 
and a similar controller was also proposed in \cite{ZhaoDorfler2015} in the context of DC microgrids. The DAPI controller can be written as 
\begin{align}
\label{eq:DAPI_vector}
\begin{aligned}
u &= -K_PV - z \\
K \dot{z} &= V - \mathcal{L}_q z,
\end{aligned}
\end{align}
where $\mathcal{L}_q$ is the Laplacian matrix of the connected graph describing the communication topology, and $K = \mathrm{diag}\{K_1,\ldots,K_n\}$ with the constant gains $K_i>0$. 
Inserting~\eqref{eq:DAPI_vector} in \eqref{eq:hvdc_vector} yields the closed-loop system
\begin{align}
\begin{aligned}
\begin{bmatrix}
K \dot{z} \\ C \dot{V}
\end{bmatrix} 
&=
\begin{bmatrix}
- \mathcal{L}_q & I_n \\
-I_n & -(\mathcal{L}_R + K_P)
\end{bmatrix}
\begin{bmatrix}
z \\ V
\end{bmatrix}.
\end{aligned}
\label{eq:hvdc_vector_DAPI}
\tag{$\mathcal{S}_\text{DAPI}$}
\end{align}

\subsection{Performance metric}
We use an $\mathcal{H}_2$ norm metric to compare the performance of the proposed controllers for MTDC grids. Consider a general input-output stable linear MIMO system~$\mathcal{S}$,
\begin{align}
\begin{aligned}
\dot x &= Ax + B w \\
y &= Hx,
\end{aligned}
\label{eq:general_lin_sys}
\tag{$\mathcal{S}$}
\end{align}
with transfer matrix $G(s) = H(sI_n - A)^{-1} B$.  The (squared) $\mathcal{H}_2$ norm of $\mathcal{S}$  is defined as
$\norm{S}_{\mathcal{H}_2}^2 \triangleq \frac{1}{2\pi} \int_{-\infty}^{\infty} \tr(G(i\omega)^*G(i\omega)) \text{d} \omega.$
The motivation for studying performance in terms of the $\mathcal{H}_2$ norm comes from two of its useful interpretations (see also \cite{Tegling2015Price}): 
\begin{enumerate}[i)]
\item If the input $w$ is a white second order process with unit covariance (white noise), then
\begin{align*}
\norm{S}_{\mathcal{H}_2}^2 = \lim_{t\rightarrow \infty} \E \{ y^*(t) y(t) \},
\end{align*}
i.e., the (squared) $\mathcal{H}_2$ norm is the steady-state variance of the output components. 

\item If the input $w\equiv 0_n$ and the initial value $x(0) = x_0$ is a random variable with covariance $\E \{ x_0 x_0^* \} = BB^*$, then
\begin{align*}
\norm{S}_{\mathcal{H}_2}^2 = \int_{0}^\infty \E \{ y^*(t) y(t) \} \text{d} t,
\end{align*}
i.e., the $\mathcal{H}_2$ norm is the 
expected $\mathcal{L}_2$ norm of the output~$y$.
\end{enumerate}
In this paper, we shall consider the following outputs when calculating the $\mathcal{H}_2$ norm:
\begin{align}
y &= \frac{1}{\sqrt{n}} V
\label{eq:output_droop}
\end{align}
for \eqref{eq:hvdc_vector_droop} and \eqref{eq:hvdc_vector_DAPI}, and
\begin{align}
y &= \frac{1}{\sqrt{n}} \tilde{V}
\label{eq:output_slack}
\end{align}
for \eqref{eq:hvdc_vector_slack}. In other words, performance will be evaluated in terms of the mean (per node) deviation from nominal voltage in the case of i) a persistent stochastic disturbance input $w$ due to e.g. fluctuating load currents, or ii) random initial voltage perturbations that are uncorrelated across nodes. 

\section{$\mathcal{H}_2$ norm evaluation}
\label{sec:computing_H2}
In this section, we compute the $\mathcal{H}_2$ norm of the
closed-loop dynamics \eqref{eq:hvdc_vector_slack}, \eqref{eq:hvdc_vector_droop}, and \eqref{eq:hvdc_vector_DAPI}. 
In order to provide tractable closed-form expressions for the \hn norms, we impose the following assumptions:
\begin{assumption}[Uniform system parameters]
\label{ass:uniform_gains}
The capacitances  $C_i$ and the gains $K_{P_i}$ and $K_i$ are uniform across the network and given by the positive constants $c$, $k_P$, and $k$, respectively. 
\end{assumption}
\begin{assumption}[Communication graph]
\label{ass:uniform_graph}
The Laplacian matrix $\mathcal{L}_q$ of the communication graph used by the DAPI controller in~\eqref{eq:hvdc_vector_DAPI} satisfies $\mathcal{L}_q = \gamma \mathcal{L}_R$, $\gamma>0$, where $\mathcal{L}_R$ was defined in \eqref{eq:hvdc_vector}.
\end{assumption}
The following theorem summarizes the main result of this section.
\begin{theorem}[Performance evaluation]
\label{th:cost_droop_slack}
Let Assumptions~\ref{ass:uniform_gains} and \ref{ass:uniform_graph} hold. 
The performance of the slack bus controlled MTDC grid~\eqref{eq:hvdc_vector_slack} with output~\eqref{eq:output_slack} is given by 
\begin{align}
\label{eq:mainthm_slack}
\norm{\mathcal{S}_\mathrm{slack}}_{\mathcal{H}_2}^2 = \frac{c}{2n}\sum_{i=1}^{n-1} \frac{1}{\tilde{\lambda}_i }, 
\end{align}
where $\tilde{\lambda}_i$ denotes the $i^{\mbox{th}}$ eigenvalue of the reduced Laplacian $\tilde{\mathcal{L}}_R$.
The performance of the droop-controlled MTDC grid~\eqref{eq:hvdc_vector_droop} with output~\eqref{eq:output_droop} is given by 
\begin{align}
\label{eq:mainthm_droop}
\norm{\mathcal{S}_\mathrm{droop}}_{\mathcal{H}_2}^2 = \frac{c}{2n} \sum_{i=1}^n \frac{1}{(\lambda_i +k_P)}, 
\end{align}
where $\lambda_i$ denotes the $i^{\mbox{th}}$ eigenvalue of $\mathcal{L}_R$. 
 The performance of the DAPI controlled MTDC grid~\eqref{eq:hvdc_vector_DAPI} with output~\eqref{eq:output_droop} is given by
\begin{align}
\label{eq:mainthm_dapi}
\norm{\mathcal{S}_\mathrm{DAPI}}_{\mathcal{H}_2}^2  =   \frac{c}{2n}  \sum_{i=1}^n  \frac{1}{ \left(\lambda_i  +  k_P  +  \frac{c\gamma \lambda_i}{c\gamma^2\lambda_i^2  + k\gamma\lambda_i^2 + kk_P\gamma\lambda_i + k} \right)}. 
\end{align}
\end{theorem}
\vspace{-2mm}
\begin{proof}
It is readily verified that the system matrices of the systems \eqref{eq:hvdc_vector_slack}, \eqref{eq:hvdc_vector_droop}, \eqref{eq:hvdc_vector_DAPI} are Hurwitz if $c,k,k_P,\gamma>0$ (noting that the grounded Laplacian $\tilde{\mathcal{L}}_R$ is positive definite if the graph $\mathcal{G}$ is connected, see \cite{Tegling2015Price}). Therefore, their \hn norms exist and are finite. To derive the respective norms, we follow the procedure outlined in \cite{Bamieh2013Price} and perform a unitary state transformation, which in the case of \eqref{eq:hvdc_vector_droop} reads $V = U\bar{V}$. The matrix $U$ diagonalizes $\mathcal{L}_R$, i.e., $\mathcal{L}_R = U^* \Lambda  U$ with $\Lambda = \diag(\lambda_1, \dots, \lambda_n)$. Since the \hn norm is unitarily invariant, we can also define $\bar{y} = U^*y$ and $\bar{w} = U^*w$ and obtain, in the new coordinates,
\begin{equation}
\begin{aligned}
c \dot{\bar{V}} &= - (\Lambda + k_P I_n) \bar{V} + \bar{w}\\
\bar{y} & = \frac{1}{\sqrt{n}}\bar{V}.
\label{eq:hvdc_vector_droop_bar}
\end{aligned}
\end{equation}
The system \eqref{eq:hvdc_vector_droop_bar} consists of $n$ decoupled subsystems~$(\mathcal{S}^i_\text{droop})$: 
\begin{equation}
\begin{aligned}
\dot{\bar{V}}_i &= - \frac{1}{c}(\lambda + k_P ) \bar{V}_i + \frac{1}{c}\bar{w} ~~\triangleq \bar{A}_\text{droop}^i \bar{V}_i + \bar{B}^i \bar{w}_i \\
\bar{y}_i & = \frac{1}{\sqrt{n}}\bar{V}_i ~~\triangleq \bar{H}^i \bar{V}_i,
\label{eq:hvdc_vector_droop_bar_i}
\end{aligned}
\end{equation}
and it holds that $\| \mathcal{S}_\text{droop}\|_{\mathcal{H}_2}^2 = \sum_{i=1}^n \|\mathcal{S}^i_\text{droop}\|_{\mathcal{H}_2}^2$. Each individual subsystem norm can now be evaluated by solving the Lyapunov equation for $P^i$:
$\bar{A}^{i*}_\text{droop} P^i + P^i \bar{A}^{i}_\text{droop}= - \bar{H}^{i*} \bar{H}^i,$
and taking $\|\mathcal{S}^i_\text{droop}\|_{\mathcal{H}_2}^2 = \mathrm{tr}(\bar{B}^{i*} P^i \bar{B}^i)$, which in our case gives
\begin{equation}
\norm{\mathcal{S}^i_\text{droop}}^2_{\mathcal{H}_2} = \frac{c}{2n(\lambda_i + k_P)}.
\end{equation}
Summing over the $n$ subsystems yields the result in \eqref{eq:mainthm_droop}. 

The $\mathcal{H}_2$ norms of \eqref{eq:hvdc_vector_slack} and  \eqref{eq:hvdc_vector_DAPI} are calculated in a similar manner. 
\end{proof}

\vspace{-2mm}
Next, we show that the $\mathcal{H}_2$ norm of
\eqref{eq:hvdc_vector_DAPI} is strictly smaller than that of \eqref{eq:hvdc_vector_droop}, which in turn is strictly smaller than that of \eqref{eq:hvdc_vector_slack}. 
\vspace{-1.8mm}
\begin{corollary}
\label{cor:slack_larger_droop}
For any choice of the parameters $c,k,k_P, \gamma>0$, it holds that
\begin{align*}
\norm{\mathcal{S}_\mathrm{DAPI}}_{\mathcal{H}_2}^2 <
\norm{\mathcal{S}_\mathrm{droop}}_{\mathcal{H}_2}^2 < \norm{\mathcal{S}_\mathrm{slack}}_{\mathcal{H}_2}^2.
\end{align*}
\end{corollary}
\vspace{-3mm}
\begin{proof}
By Cauchy's interlacing theorem 
\cite{haemers1995interlacing}, the eigenvalues $\lambda_i$ of $\mathcal{L}_R$ and the eigenvalues $\tilde{\lambda}_i$ of $\tilde{\mathcal{L}}_R$ satisfy
\begin{align*}
0 = \lambda_1 < \tilde{\lambda}_1 \le \lambda_2 \le \dots \le  \tilde{\lambda}_{n-1} \le \lambda_n.
\end{align*}
Thus
\begin{align}
\label{eq:interlacing}
\sum_{i=1}^{n-1} \frac{c}{2n\tilde{\lambda}_i } \ge \sum_{i=2}^{n} \frac{c}{2n{\lambda}_i } > \sum_{i=1}^n \frac{c}{2n(\lambda_i +k_P)}, 
\end{align}
which proves the second inequality. Furthermore, since $c,k,k_p, \gamma,\lambda_i>0$, each term of the sum in \eqref{eq:mainthm_dapi} is smaller than the corresponding term in \eqref{eq:mainthm_droop}, so $\|\mathcal{S}_\text{DAPI}\|_{\mathcal{H}_2}^2 < \|\mathcal{S}_\text{droop}\|_{\mathcal{H}_2}^2$, which concludes the proof. 
\end{proof}

\section{Performance scaling in lattice networks}
\label{sec:large-scale}
By Corollary~\ref{cor:slack_larger_droop} we know that the $\mathcal{H}_2$ norms of the droop or DAPI controlled MTDC grids are always smaller than that of the slack bus controlled grid. It turns out, as we will show in this section, that this difference in performance becomes increasingly pronounced as the network size grows. In particular, for specific network topologies, we show that the \hn norm of the slack bus controlled MTDC network grows unboundedly with network size, while it remains bounded with droop or DAPI control. This scaling of performance with network size becomes a particularly important question as DC microgrids are gaining interest, since these are likely to comprise a high number of buses.

In order to derive the relevant performance bounds, we first make the following physically motivated assumption:

\vspace{-1mm}
\begin{assumption}[Uniformly bounded resistances]
\label{ass:resistance_bounds}
The network line resistances are uniformly bounded as
\begin{align*}
\underline{R} \le R_{ij} \le \overline{R}, ~~  (i,j)\in \mathcal{E},
\end{align*}
where $\underline{R} $ and $\overline{R}$ are positive constants.  
\end{assumption}

Recall that, by \eqref{eq:mainthm_slack} and \eqref{eq:interlacing}, the $\mathcal{H}_2$ norm of the slack bus controlled MTDC grid satisfies
\begin{align}
\label{eq:slacknorm_restated}
\norm{\mathcal{S}_\text{slack}}_{\mathcal{H}_2}^2 = \frac{c}{2n}\sum_{i=1}^{n-1} \frac{1}{\tilde{\lambda}_i } \ge \frac{c}{2n}\sum_{i =2}^n \frac{1}{\lambda_i}.
\end{align}
We notice immediately that this expression blows up if one or more of the Laplacian eigenvalues $\lambda_i$ approaches zero. This is typically the case when networks grow large, unless the network is well-interconnected.  More precisely, the \hn norm's scaling will depend on how 
\begin{align}
\label{eq:defKstar}
K^*& \triangleq \frac{1}{n} \sum_{i=2}^n \frac{1}{\lambda_i}
\end{align}
scales with network size. The quantity $K^*$ is closely related to the \textit{Kirchhoff index}, also called total effective resistance or Wiener index, in resistor networks. Namely, given a network of resistors, define the pairwise effective resistance between two nodes $i$ and $j$ as $R^{\mathrm{eff}}_{ij}$. The Kirchhoff index is defined as
\begin{align}
\label{eq:defKirchhoff}
K_f &\triangleq \sum_{i<j} R^{\mathrm{eff}}_{ij},
\end{align}
see, e.g., \cite{ Lukovits1999, GhoshBoyd2008}.  
This has been shown in~\cite{Gutman1996} to equal:
\begin{align*}
K_f &= n \sum_{i=2}^n \frac{1}{\lambda_i}.
\end{align*}
Clearly $K^* =  \frac{K_f}{n^2}$, so $K^* \rightarrow \infty$ if and only if $\frac{K_f}{n^2} \rightarrow \infty$ as $n\rightarrow \infty$. 
The following well-established result proves to be very useful for our analysis. 

\vspace{-1mm}
\begin{lemma}[Rayleigh’s monotonicity law]
\label{lemma:Rayleigh}
Removing an edge from a graph, or increasing its resistance, can only increase the effective resistance between any two points in the network. Conversely, adding edges or decreasing their resistance can only decrease the effective resistance between any two points. 
\end{lemma} \vspace{-3.8mm}
\begin{proof}
See, e.g., \cite{Doyle2000}. 
\end{proof} 
\vspace{-1mm}
Rayleigh’s monotonicity law implies that well-interconnected networks have a lower Kirchhoff index, and hence a lower $\mathcal{H}_2$ norm (better performance) than sparsely interconnected networks.  It also implies that the Kirchhoff index of any network that can be \textit{embedded} in a larger network (that is a subgraph of the larger network) is lower bounded by that of the larger network~\cite{Barooah}.

We will consider a subclass of graphs for which asymptotic (in network size) bounds on the Kirchhoff index~$K_f$ can be obtained analytically, namely infinite lattices and their fuzzes, which are defined below. 


\vspace{-1mm}
\begin{definition}[Lattice]
A $d$-dimensional lattice is a graph that has a node at every point
in $\mathbb{Z}^d$ and an edge between any two nodes between which the Euclidean distance is $1$. 
\end{definition}\vspace{-3mm}

\begin{definition}[$h$-fuzz]
The $h$-fuzz of a lattice is obtained by adding an edge between any nodes within graph distance $h$.
\end{definition}
\vspace{-1mm}
\noindent By the reasoning above, the Kirchhoff index, and thereby the performance scaling, that we derive for lattices and fuzzes will provide a lower bound for all graphs that can be embedded in them. In this context, it is useful to think of the lattice dimension $d$ a measure of the graph's connectivity, which determines how performance scales in the network. 
Consider the following theorem:
\vspace{-1mm}
\begin{theorem}[Asymptotic performance in lattices]
\label{th:Kirchhoff_lattice}
Let the graph $\mathcal{G}$ corresponding to the MTDC network be a lattice or its $h$-fuzz in $d$ dimensions. Then, under Assumptions~\ref{ass:uniform_gains}--\ref{ass:resistance_bounds}, the asymptotic scaling of the \hn norm of the slack bus controlled MTDC network \eqref{eq:hvdc_vector_slack} is given in Table~\ref{table:asymptotic_performance}. This implies that for 1- and 2-dimensional lattices and $h$-fuzzes, the \hn norm scales as $n$ and $\log(n)$, respectively, and thus grows unboundedly as $n \rightarrow \infty$.

The $\mathcal{H}_2$ norms of droop and DAPI controlled MTDC networks are, on the other hand, upper bounded by 
\begin{align*}
\norm{S_\mathrm{droop}}_{\mathcal{H}_2}^2, \norm{S_\mathrm{DAPI}}_{\mathcal{H}_2}^2 &\le \frac{c}{2k_P}
\end{align*}
for any underlying network structure (that is, not limited to lattices and fuzzes). 

\end{theorem}
\vspace{-4mm}
\begin{proof}
In order to bound the \hn norm of the slack bus controlled system, by \eqref{eq:slacknorm_restated} it suffices to bound the quantity $\sum_{i=2}^{n} \frac{c}{2n{\lambda}_i } = \frac{1}{2n^2}K_f$. Note that, by  Assumption~\ref{ass:resistance_bounds}, the network's resistances are upper and lower bounded. Consider therefore the graphs where the resistances $R_{ij}$ are replaced by their lower and upper bounds. By Lemma~\ref{lemma:Rayleigh}, the Kirchhoff indices of those graphs bound the Kirchhoff index of the original graph. 
Now, by \cite{Barooah} the effective resistance between any two points $i$ and $j$ can be bounded as
\begin{align*}
\alpha_1 d_G(i,j)&\le R^{\mathrm{eff}}_{ij} \le \beta_1 d_G(i,j) \\
\alpha_2 \log(d_G(i,j))&\le  R^{\mathrm{eff}}_{ij} \le \beta_2 \log(d_G(i,j)) \\
\alpha_3 &\le  R^{\mathrm{eff}}_{ij} \le \beta_3,
\end{align*}
\noindent for, respectively, the $1$, $2$ and $3$-dimensional lattice or fuzz with uniform resistances. Here, the $\alpha$'s and $\beta$'s are positive constants, which depend on $\underline{R}$ and $\overline{R}$, but not on~$n$. The function~$d_G(i,j)$ denotes the graph distance between nodes $i$ and~$j$, which is equal to the $\ell_1$-norm between $i$ and~$j$. 

For $d=1$, the graph distance between two arbitrary nodes in a lattice with $n$ nodes, is proportional to $n$ \cite{Barooah}. Summing over all $i\le j$ as in~\eqref{eq:defKirchhoff} yields $\alpha'_1 n^3 \le K_f \le \beta'_1 n^3 \Leftrightarrow \alpha'_1 n \le K^* \le  \beta'_1 n$, for some constants $\alpha'_1,~\beta'_1$. Based on~\eqref{eq:slacknorm_restated}, the \hn norm can then be lower bounded as in Table~\ref{table:asymptotic_performance}. 

For $d=2$, the graph distance between two arbitrary nodes in a lattice with $n$ nodes, scales as $\sqrt{n}$. Summing over $i\le j$ yields $\alpha'_2 n^2 \log(n) \le K_f \le \beta'_2 n^2 \log(n)\Leftrightarrow \alpha'_2 \log(n) \le K^*\le \beta'_2 \log(n)$ for some $\alpha'_2,~\beta'_2$. The lower bound of the \hn norm can then be stated as in Table~\ref{table:asymptotic_performance}. 

For $d\ge3$, $R_{ij}$ is bounded by positive constants, so summing over $i\le j$ yields $\alpha'_3 n^2 \le K_f \le \beta'_3 n^2 \Leftrightarrow \alpha'_3 \le K^* \le \beta'_3$, that is, the \hn norm is bounded. 

Next, it is straightforward to show that $\| \mathcal{S}_\text{droop}\|^2_{\mathcal{H}_2}$ and $\|\mathcal{S}_\text{DAPI}\|^2_{\mathcal{H}_2}$ are upper bounded by $\frac{c}{2k_P}$, regardless of the $\lambda_i$, that is, of the network structure. 
\end{proof}

\begin{table}
\caption{Performance scaling of the slack bus controlled MTDC network \eqref{eq:hvdc_vector_slack} on lattices and $h$-fuzzes, as the network size $n\rightarrow \infty$, where $a_1,a_2$ are positive constants.}
\begin{center}
\renewcommand{\arraystretch}{2.0}%
\begin{tabular}{|c|c|}
\hline
{\bf Lattice dimension} & {\bf $\mathcal{H}_2$ norm, slack bus}   \\
\hline
$d =1$ & $a_1 n \le \norm{S_\text{slack}}_{\mathcal{H}_2}^2 $  \\
\hline
$d=2$ & $a_2 \log(n) \le \norm{S_\text{slack}}_{\mathcal{H}_2}^2 $   \\
\hline
\end{tabular}
\end{center}
\label{table:asymptotic_performance}
\end{table}

\begin{figure*}[h!]
  \centering
  \subfloat[][Slack bus control, $n = 10$]{
	\begin{tikzpicture}[thick,scale=0.8, every node/.style={scale=0.8}]
	\begin{axis}
	[cycle list name=Power_frequency,
	xlabel={$t$ [s]},
	ylabel={$V(t)$ [V]},
	xmin=0,
	xmax=30,
	ymin=-1,
	ymax=1,
	yticklabel style={/pgf/number format/.cd,
		fixed,
		precision=4},
	grid=major,
	height=3.8cm,
	width=0.4\textwidth,
	]
	\foreach \x in {1, 2,...,10}{
	\addplot table[x index=0,y index=\x,col sep=space]{Simulations/Plot_data/slack_string_small_n_V.txt};
	}	
	\end{axis}
	\end{tikzpicture}	
  \label{fig:a}  }
  \subfloat[][Droop control, $n = 10$] {
	\begin{tikzpicture}[thick,scale=0.8, every node/.style={scale=0.8}]
	\begin{axis}
	[cycle list name=Power_frequency,
	xlabel={$t$ [s]},
	ylabel={$V(t)$ [V]},
	xmin=0,
	xmax=30,
	ymin=-1,
	ymax=1,
	yticklabel style={/pgf/number format/.cd,
		fixed,
		precision=4},
	grid=major,
	height=3.8cm,
	width=0.4\textwidth,
	]
	\foreach \x in {1, 2,...,10}{
	\addplot table[x index=0,y index=\x,col sep=space]{Simulations/Plot_data/droop_string_small_n_V.txt};
	}	
	\end{axis}
	\end{tikzpicture}	
  \label{fig:b}}
   \subfloat[][DAPI control, $n = 10$] {
  \begin{tikzpicture}[thick,scale=0.8, every node/.style={scale=0.8}]
	\begin{axis}
	[cycle list name=Power_frequency,
	xlabel={$t$ [s]},
	ylabel={$V(t)$ [V]},
	xmin=0,
	xmax=30,
	ymin=-1,
	ymax=1,
	yticklabel style={/pgf/number format/.cd,
		fixed,
		precision=4},
	grid=major,
	height=3.8cm,
	width=0.4\textwidth,
	]
	\foreach \x in {1, 2,...,10}{
	\addplot table[x index=0,y index=\x,col sep=space]{Simulations/Plot_data/DAPI_string_small_n_V.txt};
	}	
	\end{axis}
	\end{tikzpicture}	
  \label{fig:c}} \\
    \subfloat[][Slack bus control, $n = 100$]{
  \begin{tikzpicture}[thick,scale=0.8, every node/.style={scale=0.8}]
	\begin{axis}
	[
	cycle list name=Power_frequency,
	xlabel={$t$ [s]},
	ylabel={$V(t)$ [V]},
	xmin=0,
	xmax=30,
	ymin=-1,
	ymax=1,
	yticklabel style={/pgf/number format/.cd,
		fixed,
		precision=4},
	grid=major,
	height=3.8cm,
	width=0.4\textwidth,
]
	\foreach \x in {1, 2,...,10}{
	\addplot table[x index=0,y index=\x,col sep=space]{Simulations/Plot_data/slack_string_large_n_V.txt};
	}	
	\end{axis}
	\end{tikzpicture}	
  \label{fig:d}  }
  \subfloat[][Droop control, $n = 100$] {
 	\begin{tikzpicture}[thick,scale=0.8, every node/.style={scale=0.8}]
	\begin{axis}
	[
	cycle list name=Power_frequency,
	xlabel={$t$ [s]},
	ylabel={$V(t)$ [V]},
	xmin=0,
	xmax=30,
	ymin=-1,
	ymax=1,
	yticklabel style={/pgf/number format/.cd,
		fixed,
		precision=4},
	grid=major,
	height=3.8cm,
	width=0.4\textwidth,
]
	\foreach \x in {1, 2,...,10}{
	\addplot table[x index=0,y index=\x,col sep=space]{Simulations/Plot_data/droop_string_large_n_V.txt};
	}	
	\end{axis}
	\end{tikzpicture}	
  \label{fig:e}}
   \subfloat[][DAPI control, $n = 100$] {
	\begin{tikzpicture}[thick,scale=0.8, every node/.style={scale=0.8}]
	\begin{axis}
	[
	cycle list name=Power_frequency,
	xlabel={$t$ [s]},
	ylabel={$V(t)$ [V]},
	xmin=0,
	xmax=30,
	ymin=-1,
	ymax=1,
	yticklabel style={/pgf/number format/.cd,
		fixed,
		precision=4},
	grid=major,
	height=3.8cm,
	width=0.4\textwidth,
]
	\foreach \x in {1, 2,...,10}{
	\addplot table[x index=0,y index=\x,col sep=space]{Simulations/Plot_data/DAPI_string_large_n_V.txt};
	}	
	\end{axis}
	\end{tikzpicture}	
  \label{fig:f}}
      \caption{ Voltage trajectories for a subset of the buses in a radial MTDC network (1-dimensional lattice) of size $n = 10$ and $n = 100$ with, respectively, slack bus, droop, and DAPI control. In the slack bus case, performance is substantially worse for the larger network~(d) compared to the smaller (a), which is not in the case with droop or DAPI control; compare (b--c) to (e--f). The oscillatory behavior in the DAPI case is a typical feature of integral control.  }
\label{fig:simulation}
\end{figure*}
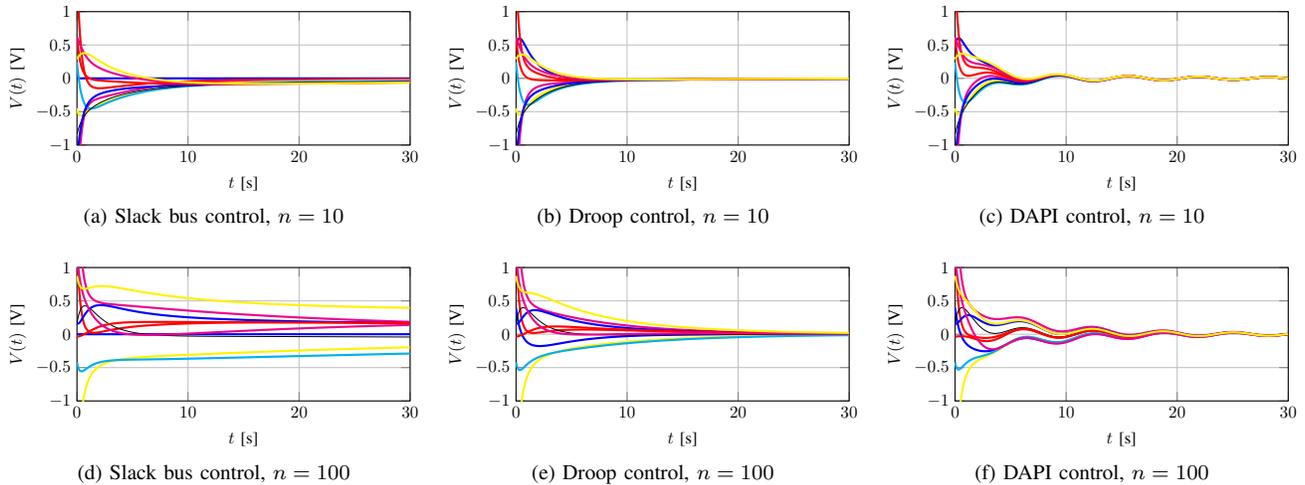

\vspace{-1mm}
Theorem~\ref{th:Kirchhoff_lattice} implies that, unless the network is well-interconnected, controlling the voltage of large MTDC grids by means of a slack bus may be unsuitable, as the \hn norm may scale unboundedly with the network size. As a result, one would obtain growing expected deviations from desired voltage levels, and thus fail in the control objective. The $\mathcal{H}_2$ norms of the droop or DAPI controlled MTDC grids, on the other hand, are uniformly upper bounded with respect to the network size, regardless of the network topology. This makes droop and DAPI control more \textit{scalable} as control strategies. 

Similar analyses as the above have been carried out for AC power grids and coupled oscillator networks in \cite{Grunberg2016, SiamiMotee2016, Tegling2016_ACC}. However, our results for the MTDC grids partly differ from those derived for AC grids, as the droop controller alone suffices to uniformly bound the \hn norm, whereas a DAPI controller is required in AC grids. This can be understood by studying the control law \eqref{eq:hvdc_vector_droop} and noting that the diagonal matrix $K_P$ provides \textit{absolute} feedback from the voltage deviations. Such absolute feedback has been shown in \cite{bamieh2012coherence} to be key in achieving bounded \hn norm scalings in systems of this type. In AC networks, absolute feedback from phase angles is not present in the droop control law, but can instead be emulated by the DAPI controller \cite{Tegling2016_ACC}.  



\section{Simulations}
\label{sec:simulations}
In this section, we demonstrate the implications of the performance scaling in Theorem~\ref{th:Kirchhoff_lattice} by a simulation of a radial MTDC grid (that is, a topology corresponding to a $1$-dimensional lattice) with sizes $n = 10$ and $n = 100$. For simplicity, we let the resistances of all DC lines be $1$~$\Omega$, and the capacitances of the DC buses be $1$~mF. The controller parameters were set to $k_P=0.1$, $k=100$ and $\gamma=1000$, respectively. The system was initiated at the voltage $V(0)=V_0$, where $V_0~\sim \mathcal{N}(0,1)$. 

The responses of the MTDC grid with the different controllers are shown in Fig.~\ref{fig:simulation}. From the figures, it is evident that the performance of the MTDC grid with a slack bus deteriorates significantly with increasing network size (greater inter-nodal differences, longer transient). The performance of the droop and DAPI controlled MTDC grids, on the other hand, remains almost unchanged as the network size increases. 
\section{Summary and Conclusions}
\label{sec:discussion}
We characterized transient performance in MTDC networks using the $\mathcal{H}_2$ norm, which can, for example, be interpreted as the expected $\mathcal{L}_2$ norm of nodal voltage deviations under Gaussian initial conditions. The performance of an MTDC grid controlled by means of a slack bus whose voltage is kept constant, is shown to be worse than that of a droop controlled grid. We showed that performance can be further improved by a DAPI controller. For network topologies resembling 1- or 2-dimensional lattices, the $\mathcal{H}_2$ norm scales unfavorably and grows unboundedly with the number of buses. On the other hand, the $\mathcal{H}_2$ norms of MTDC grids controlled with droop or DAPI controllers are always uniformly bounded with respect to the network size, regardless of topology. These control laws are therefore scalable, and thus amenable to larger MTDC grids.

\bibliographystyle{ieeetr}
\bibliography{references}
\vspace{1mm}
\end{document}